\documentclass[letterpaper, 10 pt, conference]{ieeeconf}  

\IEEEoverridecommandlockouts                              
\usepackage{graphicx} 

\usepackage{amsfonts,amsmath,amssymb} 

\usepackage{array}
\usepackage{subfigure}

\graphicspath{{Figures/}}

\usepackage[utf8]{inputenc}
\usepackage[T1]{fontenc}
\usepackage[english]{babel} 

\usepackage{calc}

\usepackage{cite}
\usepackage{subfigure}

\usepackage{amsmath}
\usepackage{amssymb}
\usepackage{eucal}

\newcommand{\barCl}{\begin{IEEEeqnarray}{rCl}}
\newcommand{\earCl}{\end{IEEEeqnarray}}
\newcommand{\ans}{\IEEEeqnarraynumspace}
\newcommand{\ase}{\IEEEyessubnumber}
\newcommand{\ane}{\IEEEnonumber}

\newcommand{\LL}{\mathcal{L}}
\newcommand{\HH}{\mathcal{H}}

\newcommand{\I}{\mathcal{I}}
\newcommand{\J}{\mathcal{J}}
\newcommand{\K}[1]{\mathcal{K}(#1)}
\newcommand{\OO}{\mathcal{O}}
\renewcommand{\P}{\mathcal{P}}
\newcommand{\Q}{\mathcal{C}}
\newcommand{\R}[1]{\mathcal{R}(#1)}
\newcommand{\dR}[1]{\mathcal{R}'(#1)}
\renewcommand{\S}{\mathcal{S}}

\newcommand{\abc}{^{abc}}
\newcommand{\ABO}{^{\alpha\beta 0}}
\newcommand{\AB}{^{\alpha\beta}}
\newcommand{\dq}{^{dq}}
\newcommand{\dqO}{^{dq0}}

\newcommand{\us}{u_s}

\newcommand{\is}{\imath_s}
\newcommand{\ir}{\imath_r}
\newcommand{\phis}{\phi_s}
\newcommand{\phir}{\phi_r}
\newcommand{\phiM}{\phi_M}
\newcommand{\Te}{T_e}
\newcommand{\Tl}{T_L}

\newcommand{\Tderive}[1]{\frac{d#1}{dt}}
\newcommand{\Pderive}[2]{\frac{\partial#1}{\partial#2}}

\newcommand{\transpose}[1]{{#1}^T}
\makeatletter
\DeclareMathOperator{\op@cos}{cos}
\DeclareMathOperator{\op@sin}{sin}
\renewcommand{\cos}[1]{\op@cos#1}
\renewcommand{\sin}[1]{\op@sin#1}
\makeatother

\newcommand{\Rs}{R_s}
\newcommand{\Rr}{R_r}
\newcommand{\Lm}{L_m}

\newcommand{\Ls}{L_s}
\newcommand{\Lr}{L_r}
\newcommand{\Jl}{J}
\newcommand{\np}{n}

\newcommand{\nbdash}{\nobreakdash-\hspace{0pt}}

\title{\LARGE \bf Energy-based modeling of electric motors}
\author{Al Kassem Jebai, Pascal Combes, François Malrait, Philippe Martin and Pierre Rouchon
\thanks{A.-K.~Jebai is with Akka Technologies and is working as a consultant with Schneider Toshiba Inverter Europe, 27120~Pacy-sur-Eure,~France.
{\tt\footnotesize al-kassem.jebai@non.schneider-electric.com}}%
\thanks{P.~Combes, P.~Martin and P.~Rouchon are with the Centre Automatique et Systèmes, MINES ParisTech, 75006~Paris,~France
{\tt\footnotesize \{philippe.martin, pierre.rouchon\}@mines-paristech.fr}}%
\thanks{P.~Combes and F.~Malrait are with Schneider Toshiba Inverter Europe, 27120~Pacy-sur-Eure,~France
{\tt\footnotesize \{pascal.combes, francois. malrait\}@schneider-electric.com}}
}

\begin{document}

\maketitle
\thispagestyle{empty}
\pagestyle{empty}

\begin{abstract}
	We propose a new approach to model electrical machines based on energy considerations and construction symmetries of the motor. We detail the approach on the Permanent-Magnet Synchronous Motor and show that it can be extended to Synchronous Reluctance Motor and Induction Motor.
	Thanks to this approach we recover the usual models without any tedious computation. We also consider effects due to non-sinusoidal windings or saturation and provide experimental data.
\end{abstract}



\section{Introduction}
	Good models of electric motors are paramount for the design of control laws. The well-established linear sinusoidal models may be not accurate enough for some applications. That is why a lot of interest is shown in modeling non-linear and non-sinusoidal effects in electrical machines. Magnetic saturation modeling has become even more critical when considering sensorless control schemes with signal injection \cite{Bianchi2011,DeBelie2005,Reigosa2007,SergeDM2009MITo}.

	The linear sinusoidal models are usually derived by a microscopic analysis of the machine, see e.g.~\cite{Chiasson2005,Krause2002}. Based on such models, there has been some effort aiming at modeling torque ripple~\cite{BiancB2002IAITo,PetroOST2000PEITo,ZhuH1992MITo} and magnetic saturation~\cite{Levi1999ECITo,StumbSDHT2003IAITo}. One problem is that the models must respect the so-called reciprocity conditions~\cite{MelkeW1990IAITo} to be physically acceptable. An alternative way to model physical systems is to use the energy-based approach, see e.g.~\cite{Whittaker1937,Landau1982}, which was applied to electrical machines in~\cite{White1959,Ortega1998,NicklOEJ1997ACITo}. An energetic approach is used to convey the dynamic behavior of the machine. 
	
		
	In this paper we recover the usual linear sinusoidal models of most of the AC machines using a simple macroscopic approach based on energy considerations and construction symmetries. Choosing an adapted frame (which happens to be the usual $dq$ frame) allows us to get simple forms for the energy function. A nice feature of this approach is that it can easily include saturation or non-sinusoidal effects, and that the reciprocity conditions are automatically enforced. We also prove the modeling of saturation can actually be done in the fictitious frames $\alpha\beta$ or $dq$ provided the star-connection scheme is used; this fact is commonly used in practice but apparently never rigorously justified.
	
	This paper is organized as follows: in section \ref{sec:hamiltonian}, we apply the energy-based approach to a general Permanent Magnet Synchronous Motor (PMSM). Then in section \ref{sec:symmetries}, we use the construction symmetries to simplify the energy function of the PMSM. In sections \ref{sec:nonsin} and \ref{sec:saturation} we develop models for the non-sinusoidal or saturated PMSM. Finally in section \ref{sec:asynchronous} we shortly show this approach can be directly applied also to the Induction Machine (IM).
\section{Energy-based modeling of the PMSM}
\label{sec:hamiltonian}
		
	\subsection{Notations}
		When $x$ is a vector we denote its coordinates in the $uvw$ frame by $x^{uvw} := \transpose{\left(x^u, x^v, x^w\right)}$. When $f$ is a scalar function we denote its gradient by $\Pderive{f}{x^{uvw}} := \transpose{\left(\Pderive{f}{x^u}, \Pderive{f}{x^v}, \Pderive{f}{x^w}\right)}$; to be consistent when $f$ is a vector function, $\Pderive{f}{x^{uvw}}$ is the transpose of its Jacobian matrix.
	
	\subsection{A brief survey of energy-based modeling}
	\label{sec:hamiltonian:general}
	
		The evolution of a physical system exchanging energy through the external forces~$Q_i$ can be found by applying a variational principle to a function~$\LL$ --the so-called Lagrangian-- of its generalized coordinates~$\left\{q_i\right\}$ and their derivatives~$\left\{\dot{q_i}\right\}$, see e.g.~\cite{Whittaker1937,Landau1982},
		\begin{equation}
			\Tderive{}\Pderive{\LL}{\dot{q_i}} - \Pderive{\LL}{q_i} = Q_i.
			\label{eqn:hamiltonian:general}
		\end{equation}

		
		However \eqref{eqn:hamiltonian:general} is not in state form, which may be inconvenient. Such a state form with ~$p_i := \Pderive{\LL}{\dot{q_i}}$ and~$q_i$ as state variables can be obtained by considering the Hamiltonian function, also called the energy function,
		
		\begin{equation}
			\HH := \transpose{p} \dot{q} - \LL \label{eqn:hamiltonian:legendre}.
		\end{equation}
		Indeed the differential of~$\HH$ is
		\barCl
			d\HH &=& \transpose{p} d\dot{q} + \transpose{\dot{q}} dp - \transpose{\Pderive{\LL}{q}}dq - \transpose{\Pderive{\LL}{\dot{q}}}d\dot{q} \ane \\
			    &=& \transpose{\dot{q}} dp - \transpose{\Pderive{\LL}{q}} dq \ane \\
			  	&=& \transpose{\Pderive{\HH}{p}} dp + \transpose{\Pderive{\HH}{q}} dq, \label{eqn:hamiltonian:difflegendre}
		\earCl
		hence $\HH$ can be seen as a function of the generalized coordinates~$\{q_i\}$ and the generalized momenta~ $\{p_i\}$. As a consequence we find the so-called Hamiltonian equations
		\barCl
			\label{eqn:hamiltonian:stateform}
			\Tderive{p_i} &=& -\Pderive{\HH}{q_i} + Q_i \ase \label{eqn:hamiltonian:dynamic} \\
			\Tderive{q_i} &=& \Pderive{\HH}{p_i}, \ase \label{eqn:hamiltonian:speed}
		\earCl
		which are in state form.
		
	\subsection{Application to a PMSM in the $abc$ frame}
		\label{sec:hamiltonian:abc}
	
		For a PMSM with three identical windings the generalized coordinates are
		\begin{equation*}
			q = \transpose{(\theta, q_s^a, q_s^b, q_s^c)},
		\end{equation*}
		where $\theta$ is the (electrical) rotor angle and $q_s\abc$ are the electrical charges in the stator windings. Their derivatives are
		\begin{equation*}
			\dot{q} = \transpose{(\omega, \is^a, \is^b, \is^c)},
		\end{equation*}
		where $\omega$ is the (electrical) rotor velocity and $\is\abc$ are the currents in the stator windings. The power exchanges are:
		\begin{itemize}
			\item the electrical power~$\transpose{\us\abc} \is\abc$ provided to the motor by the electrical source, where $\us\abc$ is the vector of voltage drops across the windings; this power is associated with the generalized force $\us\abc$
			\item the electrical power~$-\Rs \transpose{\is\abc}\is\abc$ dissipated in the stator resistances~$\Rs$; it is associated with the generalized force $-\Rs\is\abc$
			\item the mechanical power $-\Tl\frac{\omega}{\np}$ dissipated in the load, where $\Tl$ is the load torque and $\np$ the number of pole pairs; it is associated with the generalized force $-\Tl$.
		\end{itemize}
		
		Applying \eqref{eqn:hamiltonian:general} and noting there is no storage of charges in an electrical motor, hence the Lagrangian function does not depend on $q_s\abc$, we find
		\barCl
			\label{eqn:hamiltonian:dynamicLabc}
			\Tderive{}\Pderive{\LL\abc}{\is\abc} &=& \us\abc - \Rs \is\abc \ase \\
			\Tderive{}\Pderive{\LL\abc}{\omega} - \Pderive{\LL\abc}{\theta} &=& - \frac{\Tl}{\np}. \ase
		\earCl
		We denote the Lagrangian function by $\LL\abc$ to underline it is considered as a function of the variables $\is\abc$. We then recover the usual equations of the PMSM, see e.g.~\cite{Krause2002,Chiasson2005}, by defining
		\barCl
			\phis\abc(\theta, \omega, \is\abc) &:=& \Pderive{\LL\abc}{\is\abc}(\theta, \omega, \is\abc) \label{eqn:hamiltonian:phiabc}\\
			\Te\abc(\theta, \omega, \is\abc) &:=& \np\Pderive{\LL\abc}{\theta}(\theta, \omega, \is\abc); \label{eqn:hamiltonian:teLabc}
		\earCl
		$\phis\abc$ can be identified with the stator flux and $\Te\abc$ with the electro-mechanical torque. Hence the specification of the Lagrangian function yields not only the dynamical equations but also the current-flux relation and the electro-mechanical coupling.
	
		To get a system in state form we define as in \eqref{eqn:hamiltonian:legendre} the Hamiltonian function
		\begin{equation}
			\HH\abc := \omega\Pderive{\LL\abc}{\omega} + \transpose{\is\abc}\Pderive{\LL\abc}{\is\abc} - \LL\abc.
			\label{eqn:hamiltonian:legendreabc}
		\end{equation}
		$\HH\abc$ can be seen as a function of the angle $\theta$, the rotor kinetic momentum $\rho := \Pderive{\LL\abc}{\omega}$ and the stator flux $\phis\abc := \Pderive{\LL\abc}{\is\abc}$; $\HH\abc$ of course does not depend on $q_s\abc$. By \eqref{eqn:hamiltonian:difflegendre}~and~\eqref{eqn:hamiltonian:stateform} we then find the state form
		\barCl
			\label{eqn:hamiltonian:dynamicHabc}
			\Tderive{\phis\abc} &=& \us\abc - \Rs \is\abc \ase \\
			\np\Tderive{\rho} &=& \Te\abc - \Tl, \ase
		\earCl
		with
		\barCl
			\is\abc(\theta, \rho, \phis\abc) &=& \Pderive{\HH\abc}{\phis\abc}(\theta, \rho, \phis\abc) \label{eqn:hamiltonian:iabc}\\
			\Te\abc(\theta, \rho, \phis\abc) &=& -n\Pderive{\HH\abc}{\theta}(\theta, \rho, \phis\abc). \label{eqn:hamiltonian:teHabc}
		\earCl
		
		In the next subsections we show this Hamiltonian formulation can be simplified by expressing it in the $\alpha\beta$ and $dq$ frames.
		
	\subsection{Hamiltonian formulation in the $\alpha\beta$ frame}
			\label{sec:hamiltonian:AB}
		
		\begin{figure}
			\includegraphics[width=\linewidth-2pt]{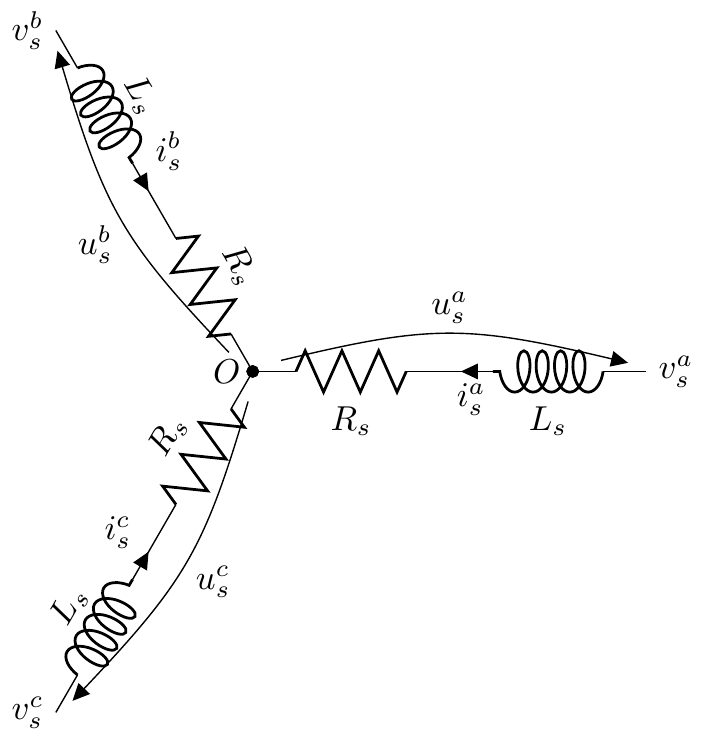}
			\caption{Star-connected motor electrical circuit}
			\label{fig:hamiltonian:star}
		\end{figure}
		
		The stator windings of the PMSMs are usually star-connected, see figure \ref{fig:hamiltonian:star}. This implies
		\begin{equation}
			\is^a + \is^b + \is^c = 0.
			\label{eqn:hamiltonian:i0}
		\end{equation}
		This algebraic relation can easily be taken into account after a change of coordinates. Indeed we change variables to the $\alpha\beta 0$ frame with $x\ABO := \Q x\abc$, thanks to the orthogonal matrix (i.e. $\Q^{-1} = \transpose{\Q}$)
		\begin{equation*}
			\Q := \sqrt{\frac{2}{3}}\begin{pmatrix}
				1 & -\frac{1}{2} & -\frac{1}{2} \\
				0 & \frac{\sqrt{3}}{2} & -\frac{\sqrt{3}}{2} \\
				\frac{1}{\sqrt{2}} & \frac{1}{\sqrt{2}} & \frac{1}{\sqrt{2}}
			\end{pmatrix}.
		\end{equation*}
		We then define the Hamiltonian function in the $\alpha\beta 0$~variables by
		\begin{equation*}
			\HH\ABO(\theta, \rho, \phis\ABO) := \HH\abc(\theta, \rho, \transpose{\Q}\phis\ABO).
		\end{equation*}
		This transformation preserves \eqref{eqn:hamiltonian:dynamicHabc}, \eqref{eqn:hamiltonian:iabc} and \eqref{eqn:hamiltonian:teHabc}; for instance
		\begin{equation*}
			\is\ABO = \Q\is\abc = \Q\Pderive{\HH\abc}{\phis\abc} = \Pderive{\HH\ABO}{\phis\ABO}
		\end{equation*}
		and
		\barCl
			\Te\ABO(\theta, \rho, \phis\ABO) &:=& \Te\abc(\theta, \rho, \Q\phis\ABO) \ane \\
			&=& -n\Pderive{\HH\abc}{\theta}(\theta, \rho, \Q\phis\ABO) \ane \\
			&=& -n\Pderive{\HH\ABO}{\theta}(\theta, \rho, \phis\ABO). \ane
		\earCl
		The constraint \eqref{eqn:hamiltonian:i0}, i.e. $\is^0(\theta, \rho, \phis\ABO) = 0$, and the assumption of a non-degenerated Hamiltonian function implies $\phis^0$ is a function of $(\theta, \rho, \phis^\alpha, \phis^\beta)$ by the implicit function theorem. Hence we can define the star-connection-constrained Hamiltonian function
		\begin{equation*}
			\HH\AB(\theta, \rho, \phis\AB) := \HH\ABO\Bigl(\theta, \rho, \bigl(\phis\AB, \phis^0(\theta, \rho, \phis\AB)\bigr)\Bigr).
		\end{equation*}
		Obviously, the system can be decomposed into
		\barCl
			\label{eqn:hamiltonian:dynamicHAB}
			\Tderive{\phis\AB} &=& \us\AB - \Rs \is\AB  \ase \\
			\np\Tderive{\rho} &=& \Te\AB - \Tl \ase \label{eqn:hamiltonian:mechanicalHAB}\\
			\ane \\
			\Tderive{\phis^0} &=& \us^0; \label{eqn:hamiltonian:dynamic0HAB}
		\earCl
		moreover
		\barCl
			\Pderive{\HH\AB}{\phis\AB} &=& \Pderive{}{\phis\AB}\HH\ABO\Bigl(\theta, \rho, \bigl(\phis\AB, \phis^0(\theta, \rho, \phis\AB)\bigr)\Bigr) \ane \\
			&=& \Pderive{\HH\ABO}{\phis\AB} + \Pderive{\HH\ABO}{\phis^0}\Pderive{\phis^0}{\phis\AB} \ane \\
			&=& \Pderive{\HH\ABO}{\phis\AB} \ane \\
			&=:& \is\AB(\theta, \rho, \phis\AB) \label{eqn:hamiltonian:iAB} \\
			-\np\Pderive{\HH\AB}{\theta} &=& -\np\Pderive{}{\theta}\HH\ABO\Bigl(\theta, \rho, \bigl(\phis\AB, \phis^0(\theta, \rho, \phis\AB)\bigr)\Bigr) \ane \\
			&=& -\np\Pderive{\HH\ABO}{\theta} - \np\Pderive{\HH\ABO}{\phis^0}\Pderive{\phis^0}{\theta} \ane \\
			&=& -\np\Pderive{\HH\ABO}{\theta} \ane \\
			&=:& \Te\AB(\theta, \rho, \phis\AB), \label{eqn:hamiltonian:TeAB}
		\earCl
		where we used $\Pderive{\HH\ABO}{\phis^0}\Bigl(\theta, \rho, \bigl(\phis\AB, \phis^0(\theta, \rho, \phis\AB)\bigr)\Bigr) = \is^0 = 0$. This means the current-flux and electromechanical relations are also decoupled from the $0$\nbdash axis.
		
		Therefore we have simplified the equation coming from the Hamiltonian formulation by decoupling from the $0$\nbdash axis (there are less equations and less variables). The derivation is valid for any Hamiltonian function, which is usually not acknowledged in the literature.
		
		
	\subsection{Hamiltonian formulation in the $dq$ frame}
		\label{sec:hamiltonian:dq}
		
		We can further simplify the formulation by expressing variables in the $dq0$ frame, i.e.~$\phis\dqO := \transpose{\R{\theta}}\phis\ABO$ with
		\begin{equation*}
			\R{\theta} := \begin{pmatrix}
				\cos{\theta} & -\sin{\theta} & 0 \\
				\sin{\theta} & \cos{\theta} & 0 \\
				0 & 0 & 1 \\
			\end{pmatrix},
		\end{equation*}
		and defining
		\begin{equation*}
			\HH\dqO(\theta, \rho, \phis\dqO) := \HH\ABO(\theta, \rho, \R{\theta}\phis\dqO).
		\end{equation*}
		Unfortunately this transformation does not preserve the Hamiltonian equations. However the flavor of the Hamiltonian formulation is preserved; indeed on the one hand
		\barCl
			\Tderive{\phis\dqO} &=& \Tderive{}\left(\transpose{\R{\theta}}\phis\ABO\right) \ane \\
							    &=& \transpose{\R{\theta}}\Tderive{\phis\ABO} + \Tderive{\transpose{\R{\theta}}}\phis\ABO \ane \\
							    &=& \transpose{\R{\theta}}(\us\ABO - \Rs \is\ABO) + \omega\transpose{\dR{\theta}}\R{\theta}\phis\dqO \ane \\
							   &=& \us\dqO - \Rs \is\dqO - \J_3 \omega \phis\dqO \ase \label{eqn:hamiltonian:dynamicHdq0} \\
			\np\Tderive{p} &=& \Te\dqO - \Tl, \ase \label{eqn:hamiltonian:mechanicalHdqO}
		\earCl
		where
		\begin{equation*}
			\J_3 := -\transpose{\dR{\theta}}\R{\theta} = \begin{pmatrix}
						0 & -1 & 0 \\
						1 & 0 & 0 \\
						0 & 0 & 0 \\
					\end{pmatrix}.
		\end{equation*}
		On the other hand
		\barCl
			\Pderive{\HH\dqO}{\phis\dqO} &=& \Pderive{\phis\ABO}{\phis\dqO} \Pderive{\HH\ABO}{\phis\ABO} \ane \\
			&=& \transpose{\R{\theta}} \is\ABO \ane \\
			&=:& \is\dqO \ane \\
			\Pderive{\HH\dqO}{\theta} &=& \Pderive{\HH\ABO}{\theta} + \transpose{\Pderive{\phis\ABO}{\theta}}\Pderive{\HH\ABO}{\phis\ABO} \ane \\
			&=& \Pderive{\HH\ABO}{\theta} + \transpose{\left(\dR{\theta}\phis\dqO\right)}  \R{\theta}\Pderive{\HH\dqO}{\phis\dqO} \ane \\
			&=& \Pderive{\HH\ABO}{\theta} - \transpose{\phis\dqO} \J_3 \is\dqO, \ane
		\earCl
		hence the current-flux relation and electro-mechanical torque are
		\barCl
			\is\dqO(\theta, \rho, \phis\dqO) &=& \Pderive{\HH\dqO}{\phis\dqO}(\theta, \rho, \phis\dqO) \label{eqn:hamiltonian:idq0} \\
			\Te\dqO(\theta, \rho, \phis\dqO) &:=& \Te\ABO(\theta, \rho, \R{\theta} \phis\dqO) \ane \\
												&=& -\np\Pderive{\HH\dqO}{\theta} + \np\transpose{\is\dqO} \J_3 \phis\dqO. \label{eqn:hamiltonian:Tedq0}	
		\earCl
		
		Since $\is^0(\theta, \rho, \phis\dqO) = 0$ when evaluated under the constraint \eqref{eqn:hamiltonian:i0}, the $0$\nbdash axis can be decoupled as in section \ref{sec:hamiltonian:AB}:
		\barCl
			\label{eqn:hamiltonian:dynamicHdq}
			\Tderive{\phis\dq} &=& \us\dq - \Rs \is\dq - \J \omega \phis\dq \ase \\
			\np\Tderive{\rho} &=& \Te\dq - \Tl \ase \label{eqn:hamiltonian:mechanicalHdq}\\
			\ane \\
			\Tderive{\phis^0} &=& \us^0, \label{eqn:hamiltonian:dynamic0Hdq}
		\earCl
		with current-flux relation and electro-mechanical torque given by
		\barCl
			\is\dq(\theta, \rho, \phis\dq) &=& \Pderive{\HH\dq}{\phis\dq}(\theta, \rho, \phis\dq) \label{eqn:hamiltonian:idq} \\
			\Te\dq(\theta, \rho, \phis\dq) &=& -\np\Pderive{\HH\dq}{\theta} + \np\transpose{\is\dq} \J \phis\dq \label{eqn:hamiltonian:Tedq}	
		\earCl
		where $\J := \begin{pmatrix}
			0 & -1 \\
			1 & 0 \\
		\end{pmatrix} $.
		
		We will see in the next section that the construction symmetries of the PMSM are more easily expressed in the $dq$ frame, resulting in simpler Hamiltonian functions.
		
	\subsection{Partial conclusion}
		The whole model of the PMSM can thus be obtained with the specification of only one energy function, yet to be defined. Since no assumption was made on the motor, this approach applies to any PMSM. In particular this implies that modeling the saturation in the $dq$ frame is equivalent to modeling it in the physical frame $abc$ if the motor is star-connected; to our knowledge this had never been proven before though the conclusion is widely used.
		
		Besides the reciprocity condition \cite{MelkeW1990IAITo} of the flux-current relation $\Pderive{\phis^d}{\is^q} = \Pderive{\phis^q}{\is^d}$ directly stems from the energy formulation. Indeed, as~$\is^d = \Pderive{\HH\dq}{\phis^d}$ and~$\is^q = \Pderive{\HH\dq}{\phis^q}$, we have
		\begin{equation*}
			\Pderive{\is^d}{\phis^q} = \frac{\partial^2\HH}{\partial\phis^q \partial\phis^d} = \frac{\partial^2\HH}{\partial\phis^d \partial\phis^q} = \Pderive{\is^q}{\phis^d},
		\end{equation*}
		which is equivalent to the reciprocity condition.

\section{Construction symmetry considerations}
	\label{sec:symmetries}

	To restrict the number of possible Hamiltonian functions we now put constraints on the form of these functions. To do so we use three simple and general geometric symmetries enjoyed by any well-built PMSM. 	
	
	
	\subsection{Phase permutation symmetry}
	\label{sec:symmetries:permutation}
		Circularly permuting the phases, then rotating the rotor by~$\frac{2\pi}{3}$ leaves the motor unchanged, hence the energy. Thus
		\barCl
			\HH\abc(\theta, \rho, \phis\abc) &=& \HH\abc(\theta+\frac{2\pi}{3}, \rho, \P\phis\abc),
			\label{eqn:symmetry:pabc}
		\earCl
		where
		\begin{equation*}
			\P := \begin{pmatrix}
				0 & 1 & 0 \\
				0 & 0 & 1 \\
				1 & 0 & 0
			\end{pmatrix}.
		\end{equation*}
		Writing this relation in the $\alpha\beta 0$ and $dq0$ frames yields
		\barCl
			\HH\ABO(\theta, \rho, \phis\ABO) &=& \HH\ABO(\theta+\frac{2\pi}{3}, \rho, \Q \P \transpose{\Q} \phis\ABO) \ans \label{eqn:symmetry:pAB} \\
			\HH\dqO(\theta, \rho, \phis^d, \phis^q, \phis^0) &=& \HH\dqO(\theta + \frac{2\pi}{3}, \rho, \phis^d, \phis^q, \phis^0).
			\label{eqn:symmetry:pdq}
		\earCl
	
	\subsection{Central symmetry}
	
		Reversing the currents in the phases, then rotating the rotor by $\pi$ leaves the motor unchanged, hence the energy. Thus
		\barCl
			\HH\abc(\theta, \rho, \phis\abc) &=& \HH\abc(\theta+\pi, \rho, -\phis\abc).
			\label{eqn:symmetry:oabc}
		\earCl
		Writing this relation in the $\alpha\beta 0$ and $dq0$ frames yields
		\barCl
			\HH\ABO(\theta, \rho, \phis\ABO) &=& \HH\ABO(\theta+\pi, \rho, -\Q \transpose{\Q} \phis\ABO)
			\label{eqn:symmetry:oAB} \ans \\
			\HH\dqO(\theta, \rho, \phis^d, \phis^q, \phis^0) &=& \HH\dqO(\theta + \pi, \rho, \phis^d, \phis^q, -\phis^0).
			\label{eqn:symmetry:odq}
		\earCl
	
	\subsection{Orientation symmetry}
		
		Permuting the phases $b$ and $c$ preserves the energy, then changing direction. the direction of rotation leaves the motor unchanged, hence the energy. Thus
		\barCl
			\HH\abc(\theta, \rho, \phis\abc) &=& \HH\abc(-\theta, -\rho, \OO\phis\abc),
			\label{eqn:symmetry:rabc}
		\earCl
		where
		\begin{equation*}
			\OO := \begin{pmatrix}
				1 & 0 & 0 \\
				0 & 0 & 1 \\
				0 & 1 & 0
			\end{pmatrix}.
		\end{equation*}
		Writing this relation in the $\alpha\beta 0$ and $dq0$ frames yields
		\barCl
			\HH\ABO(\theta, \rho, \phis\ABO) &=& \HH\ABO(-\theta, -\rho, \Q \OO\transpose{\Q} \phis\ABO) \ans
			\label{eqn:symmetry:rAB} \\
			\HH\dqO(\theta, \rho, \phis^d, \phis^q, \phis^0) &=& \HH\dqO(-\theta, -\rho, \phis^d, -\phis^q, \phis^0).
			\label{eqn:symmetry:rdq}
		\earCl
		
	\subsection{Partial conclusion}
	\label{sec:symmetries:conclusions}
		Gathering \eqref{eqn:symmetry:pdq}, \eqref{eqn:symmetry:odq} and \eqref{eqn:symmetry:rdq} and decoupling the $0$\nbdash axis, we eventually find
		\barCl
			\HH\dq(\theta, \rho, \phis^d, \phis^q) &=& \HH\dq(\theta + \frac{\pi}{3}, \rho, \phis^d, \phis^q) \ase \label{eqn:symmetry:dqperiod} \\
			\HH\dq(\theta, \rho, \phis^d, \phis^q) &=& \HH\dq(-\theta, -\rho, \phis^d, -\phis^q). \ase
			\label{eqn:symmetry:dqsymmetry}
		\earCl
		In other words, $\HH\dq$ is $\frac{\pi}{3}$\nbdash periodic with respect to $\theta$ and satisfies a parity condition on $\theta$, $\rho$ and $\phis^q$. These symmetries constrains the possible energy functions as shown in the next sections.
		
	
	\subsection{The linear sinusoidal model}
		\label{sec:symmetries:linearsinus}
	
		As an example we consider the simplest case, namely a PMSM whose magnetic energy in the $dq$ frame is a second-order polynomial not depending on the position $\theta$ nor on the kinetic momentum $\rho$. This means we assume a sinusoidally wound motor with a first-order flux-current relation. Moreover, as we are not modeling mechanics, we take the simplest kinetic energy. That is to say
		\begin{equation}
			\HH_l\dq := \frac{\rho^2}{2J\np^2}+ a + b\phis^d + c\phis^q + \frac{d}{2}{\phis^d}^2 + e\phis^d\phis^q + \frac{f}{2}{\phis^q}^2,
			\label{eqn:symmetry:quadForm}
		\end{equation}
		where J is the rotor inertia moment and $a, b, c, d, e, f$ are some constants.
		
		The symmetry \eqref{eqn:symmetry:dqsymmetry} implies $c = e = 0$. As the the energy function $\HH\dq$ is defined up to a constant we can freely change~$a$, in particular set $a = \frac{b^2}{2}$. Defining
		\begin{itemize}
			\item the $d$\nbdash axis inductance $L^d := \frac{1}{d}$
			\item the $q$\nbdash axis inductance $L^q := \frac{1}{f}$
			\item the permanent magnet flux $\phiM := L^d b$,
		\end{itemize}
		\eqref{eqn:symmetry:quadForm} eventually reads
		\begin{equation}
			\HH_l\dq = \frac{1}{2J\np^2}\rho^2 +	\frac{1}{2L^d}(\phis^d - \phi_M)^2 + \frac{1}{2L^q}{\phis^q}^2.
			\label{eqn:symmetry:linearH}
		\end{equation}
		As a consequence \eqref{eqn:hamiltonian:dynamicHdq}, \eqref{eqn:hamiltonian:idq} and \eqref{eqn:hamiltonian:Tedq} become
		\barCl
			\Tderive{\phis\dq} &=& \us\dq - \Rs \is\dq - \J \omega \phis\dq \ase \\
			\np\Tderive{\rho} &=& \Te\dq - \Tl \ase\\
			\ane
		\earCl
		\barCl
			\is^d &=& \frac{1}{L^d}(\phis^d - \phiM) \ane \\
			\is^q &=& \frac{1}{L^q}\phis^q \ane \\
			\Te\dq &=& \np\transpose{\is\dq} \J \phis\dq = \np\left(\frac{1}{L^q} - \frac{1}{L^d}\right)\phis^d\phis^q + \frac{\np}{L^d}\phis^q\phiM, \ane
		\earCl
		which is the usual model for PMSM, see e.g. \cite{Chiasson2005,Krause2002}. It is remarkable that this model can be recovered without the rather traditional microscopic approach. We have simply followed a standard energy approach with simplest possible energy function, and taken into account very general construction symmetries.
		
		Notice the model of the Synchronous Reluctance Motor can be obtained in exactly the same way. Indeed since the rotor is not oriented, we have the extra symmetry
		\begin{equation}
			\HH\dqO(\theta, \rho, \phis^d, \phis^q, \phis^0) = \HH\dqO(\theta, \rho, -\phis^d, -\phis^q, -\phis^0),
		\end{equation}
		which implies~$b = 0$ in \eqref{eqn:symmetry:quadForm} hence~$\phiM = 0$.

\section{A non-sinusoidal PMSM model}
\label{sec:nonsin}

	One interest of the energy approach is to provide models more general than the usual sinusoidal and saturated PMSM, simply by considering more general energy functions. In particular it easily explains the so-called torque ripple phenomenon, i.e. the $\frac{\pi}{3}$\nbdash periodicity of the torque with respect to $\theta$, see e.g. \cite{PetroOST2000PEITo,BiancB2002IAITo}. We still assume the magnetic energy does not depend on the kinetic momentum $\rho$, and the simplest possible kinetic energy.

	By \eqref{eqn:symmetry:dqperiod} $\HH\dq$ is $\frac{\pi}{3}$\nbdash periodic with respect to $\theta$ hence can be expended in Fourier series
	\begin{multline}
		\HH\dq(\theta, \rho, \phis^d, \phis^q) =  \frac{1}{2J\np^2}\rho^2 + \HH\dq_0(\phis^d, \phis^q) \\
		+ \sum\limits_{k=1}^{\infty} \underbrace{a_{6k}(\phis^d, \phis^q)\cos{6k\theta} + b_{6k}(\phis^d, \phis^q)\sin{6k\theta}}_{\HH\dq_{6k}}.
	\end{multline}
	
	Thanks to symmetry \eqref{eqn:symmetry:rdq} $\HH_0\dq$ and~$\left\{a_{6k}\right\}$ are even functions of~$\phis^q$, and~$\left\{b_{6k}\right\}$ are odd functions of $\phis^q$. Particularizing \eqref{eqn:hamiltonian:idq}\nbdash \eqref{eqn:hamiltonian:Tedq} to this energy function gives
	\barCl
		\is\dq(\theta, \rho, \phis) &=& \Pderive{\HH\dq_0}{\phis}(\rho, \phis\dq) + \sum\limits_{k=1}^{\infty} \Pderive{\HH\dq_{6k}}{\phis}(\theta, \rho, \phis\dq) \ane \\
		\Te\dq(\theta, \rho, \phis) &=& -\np\sum\limits_{k=1}^{\infty} \Pderive{\HH\dq_{6k}}{\theta}(\theta, \rho, \phis\dq) + n\transpose{\is\dq}\J\phis, \ane
	\earCl
	which shows $\is\dq$ and $\Te\dq$ are also $\frac{\pi}{3}$\nbdash periodic.
	
\begin{figure}
		\subfigure[SPM torque measurement]{\includegraphics[width=\linewidth/2-8pt]{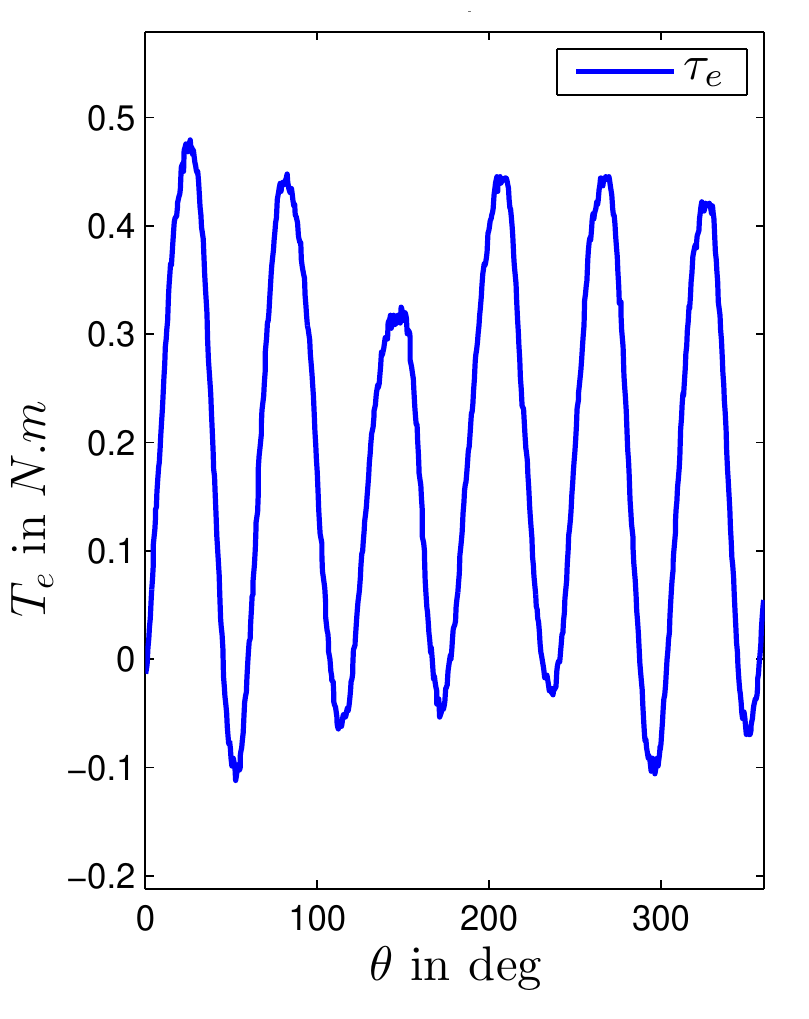} \label{fig:nonsin:6omegas:TeSPM}}
		\subfigure[IPM torque measurement]{\includegraphics[width=\linewidth/2-8pt]{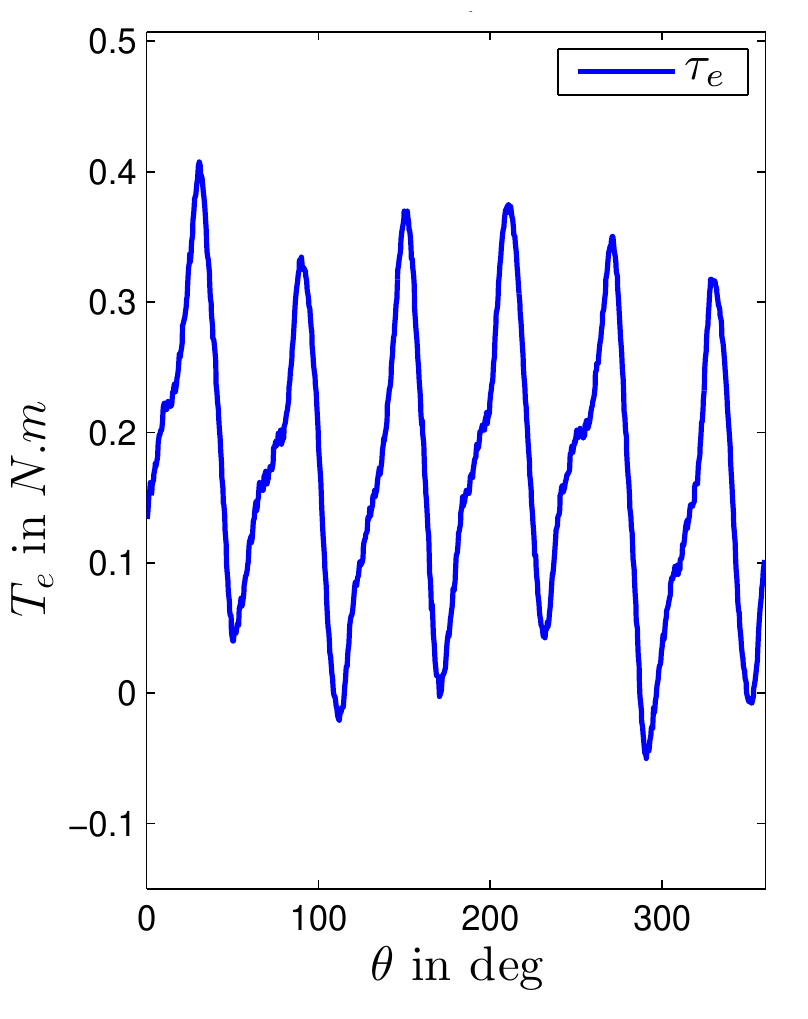} \label{fig:nonsin:6omegas:TeIPM}} \\
		\subfigure[SPM current $\is^q$ measurement]{\includegraphics[width=\linewidth/2-8pt]{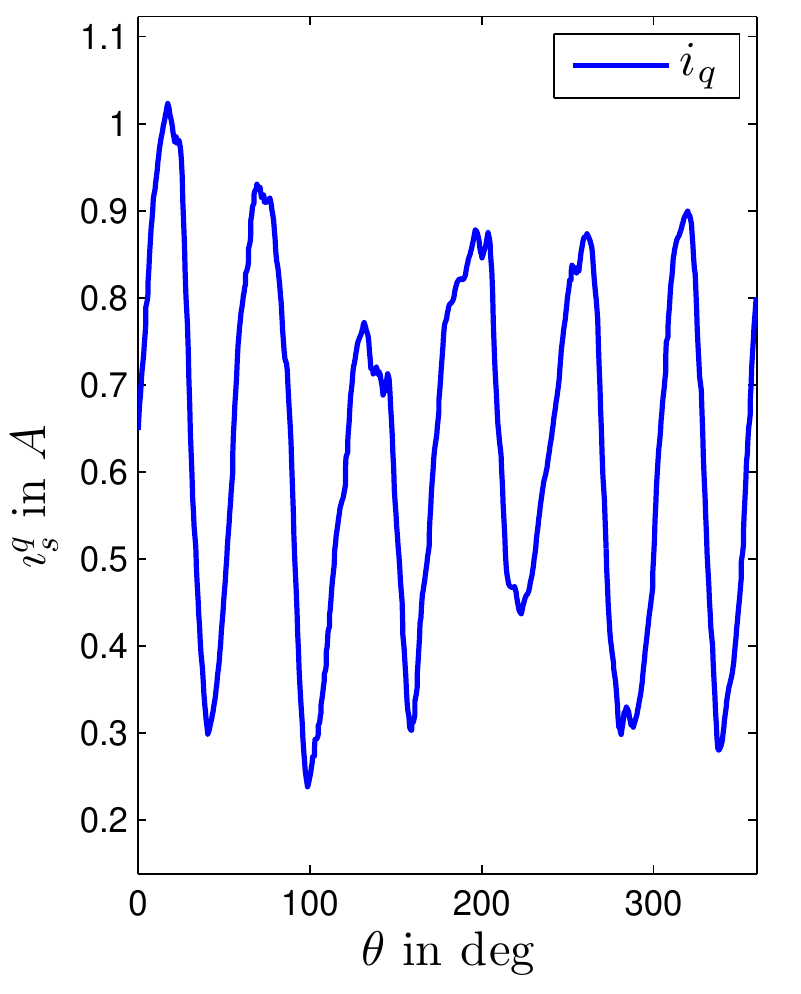} \label{fig:nonsin:6omegas:iSPM}}
		\subfigure[IPM current $\is^q$ measurement]{\includegraphics[width=\linewidth/2-8pt]{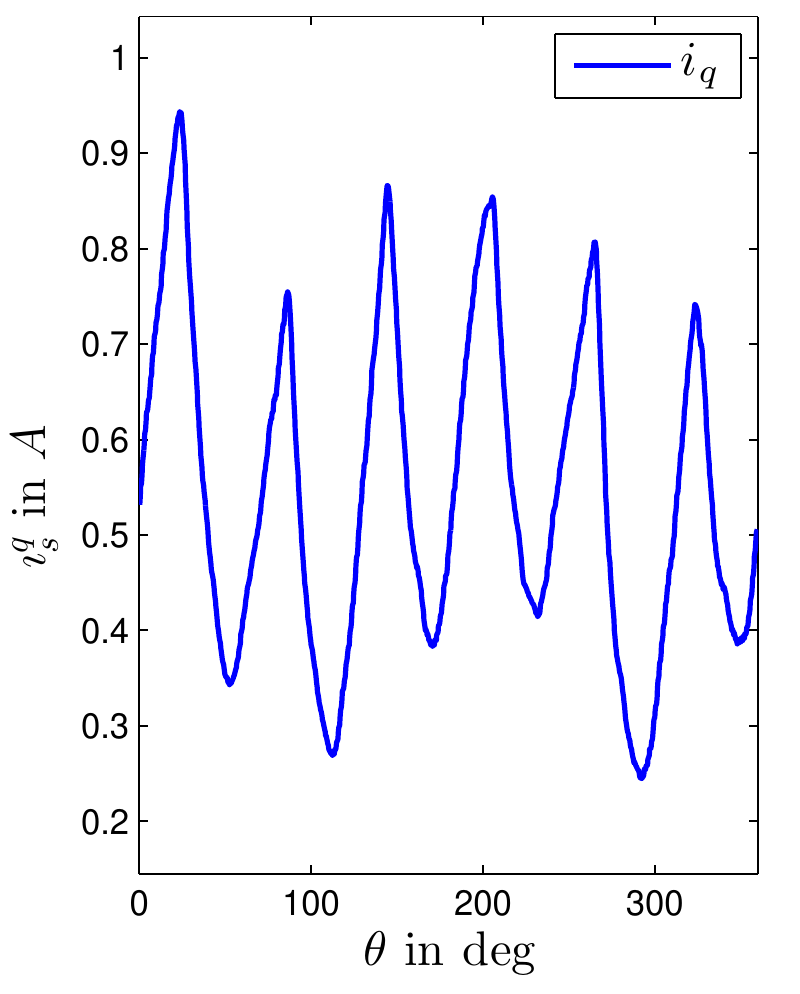} \label{fig:nonsin:6omegas:iIPM}}
		\caption{Stator current and torque measurements diverse kinds of PMSM}
		\label{fig:nonsin:6omegas}
	\end{figure}

	We experimentally checked this phenomenon on a test bench featuring current, position and torque sensors. We used two test motors, a Surface Permanent Magnet (SPM) and an Interior Permanent Magnet (IPM) PMSM, see characteristics in table \ref{tbl:nonsin:motors}. As expected the experimental plots in figure \ref{fig:nonsin:6omegas} exhibit a $\frac{\pi}{3}$\nbdash periodicity with respect to~$\theta$. The experiments  were done at low velocity and no load so that this effect is well-visible.
	

	\begin{table}
		\centering
		\renewcommand{\arraystretch}{1.25}
		\begin{tabular}{|l|l|l|}
			\hline
				PMSM kind & IPM & SPM \\
			\hline
				Rated power & $750W$ & $1500W$ \\
				Rated current (peak) & $4.51A$ & $5.19A$ \\
				Rated voltage (peak) & $110V$ & $245V$ \\
				Rotor flux (peak) & $196mWb$ & $155mWb$ \\
				Rated speed & $1800rpm$ & $3000rpm$ \\
				Rated torque & $3.98Nm$ & $6.06Nm$ \\
				Pole number ($\np$) & $3$ & $5$ \\
			\hline
		\end{tabular}
		\caption{Test motor parameters.}
		\label{tbl:nonsin:motors}
	\end{table}

	Moreover if we consider the $0$\nbdash axis, the symmetries \ref{sec:symmetries:permutation} implies $\HH\dqO$ hence $\phis^0$ is only $\frac{2\pi}{3}$\nbdash periodic with respect to~$\theta$. This effect can be experimentally seen on the potential $v_N$ of the point $O$ in figure \ref{fig:hamiltonian:star}, thanks to \eqref{eqn:hamiltonian:dynamic0Hdq}
	\begin{equation*}
		\Tderive{\phis^0}(\theta, \rho, \phis\dq) = \us^0 = v_s^0 - \sqrt{3}v_N;
	\end{equation*}
	here $v_s^0 := \frac{1}{\sqrt{3}}(v_s^a + v_s^b + v_s^c)$ is as usual set to~$0$ by the inverter. Therefore $v_N$ will exhibit a $\frac{2\pi}{3}$\nbdash periodicity with respect to~$\theta$, which was also measured on the test bench.
	
	

\section{Modeling of magnetic saturation}
\label{sec:saturation}

	We now investigate the effect of magnetic saturations; this very important when trying to control the motor at low velocity and high load, see e.g. \cite{Bianchi2011,DeBelie2005,Reigosa2007,SergeDM2009MITo}. We consider only sinusoidal motors (i.e. the energy function $\HH\dq$ is independent of $\theta$) since the non-sinusoidal effects in well-wound PMSMs are experimentally small in the presence of magnetic saturation. We still assume the magnetic energy does not depend on the kinetic momentum $\rho$, and the simplest possible kinetic energy.

	
	In normal operation $\phis^d$ is close to the permanent magnet flux $\phiM$, while $\phis^q$ is small with respect to $\phiM$. It is thus natural to expand $\HH\dq$ as a Taylor series in the variables $(\phis^d - \phiM)$ and $\phis^q$
	\begin{equation}
		\HH\dq = \HH_l\dq + \sum\limits_{n = 3}^{\infty}\sum\limits_{k = 0}^{n}\alpha_{n-k,k}(\phis^d-\phiM)^{n-k}{\phis^q}^k,
	\end{equation}
	where $\HH_l\dq$ is given by \eqref{eqn:symmetry:linearH}. Moreover, all odd powers of $\phis^q$ have by \eqref{eqn:symmetry:dqsymmetry} null coefficients, hence
	\begin{equation}
		\HH\dq = \HH_l\dq + \sum\limits_{n = 3}^{\infty}\sum\limits_{m = 0}^{\lfloor\frac{n}{2}\rfloor}\alpha_{n-2m,2m}(\phis^d-\phiM)^{n-2m}{\phis^q}^{2m}.
	\end{equation}
	
\begin{figure}
		\subfigure[IPM motor]{\includegraphics[width=\linewidth-2pt]{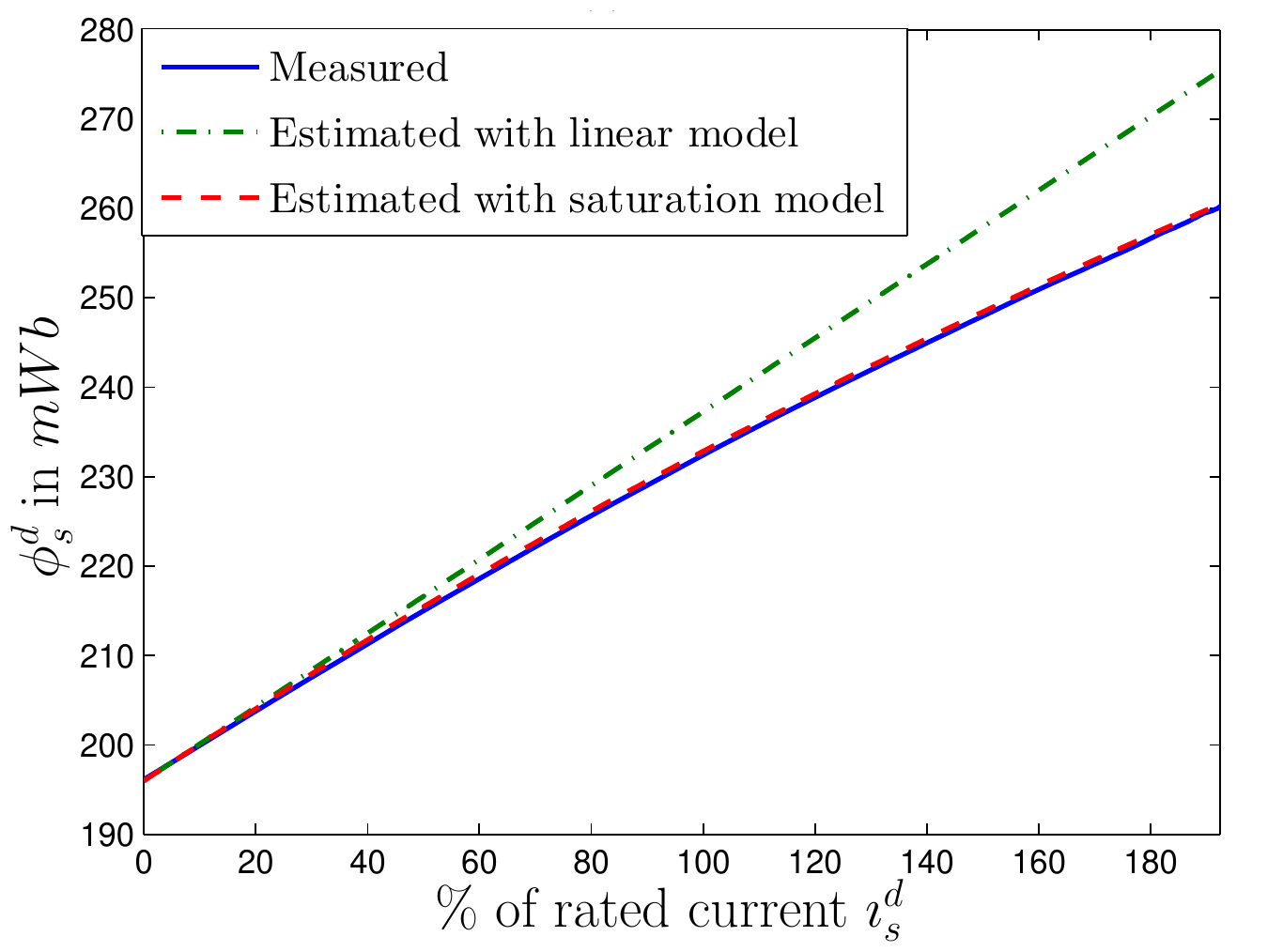}} \\
		\subfigure[SPM motor]{\includegraphics[width=\linewidth-2pt]{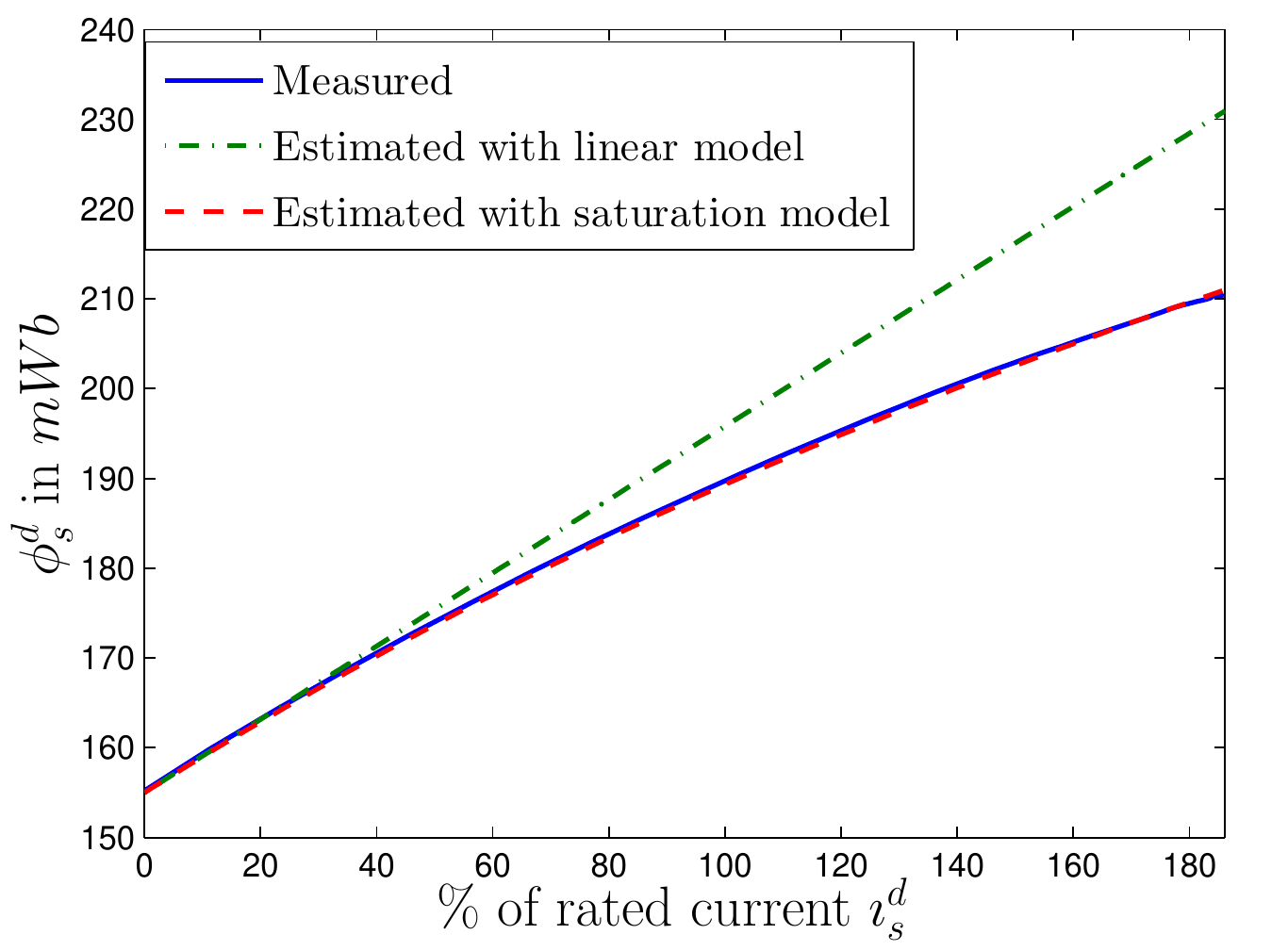}}
		\caption{Experimental and fitted flux-current relations.}
		\label{fig:saturation:measures}
	\end{figure}

	We experimentally checked the validity of this conclusion on the two motors described in table  \ref{tbl:nonsin:motors}. We first obtained the flux-current relation by integrating the back-electromotive force when applying voltage steps, see figure \ref{fig:saturation:measures}. We then truncated the series at $n = 4$ and experimentally identified $L^d, L^q, \alpha_{3,0}, \alpha_{1,2}, \alpha_{4,0}, \alpha_{2,2}, \alpha_{0,4}$, see \cite{Jebai2011} for details. The agreement between the flux-current relation obtained from $\HH\dq$ and the experimental flux-current relation is excellent. Notice the linear model using only $\HH_l\dq$ is good only at low current.
	
	\begin{table}
		\centering
		\renewcommand{\arraystretch}{1.5}
		\begin{tabular}{|l|l|l|}
			\hline
			Motor & IPM & SPM \\
			\hline
			Measured $\Rs$ & $1.52\Omega$ & $2.1\Omega$ \\
			$\frac{\phiM^2}{L^d}$ & $4.20\pm 0.12A.Wb$ & $3.06\pm 0.08A.Wb$ \\
			$\frac{\phiM^2}{L^q}$ & $2.83\pm 0.12A.Wb$ & $2.94\pm 0.08A.Wb$ \\
			$\phiM^3\alpha_{3,0}$ & $0.770\pm 0.007A.Wb$ & $0.655\pm 0.006A.Wb$ \\
			$\phiM^3\alpha_{1,2}$ & $0.702\pm 0.009A.Wb$ & $0.617\pm 0.010A.Wb$ \\
			$\phiM^4\alpha_{4,0}$ & $0.486\pm 0.012A.Wb$ & $0.724\pm 0.010A.Wb$ \\
			$\phiM^4\alpha_{2,2}$ & $0.734\pm 0.015A.Wb$ & $1.010\pm 0.025A.Wb$ \\
			$\phiM^4\alpha_{0,4}$ & $0.175\pm 0.004A.Wb$ & $0.262\pm 0.006A.Wb$ \\
			\hline
		\end{tabular}
		\caption{Experimental magnetic parameters}
		\label{tbl:saturation:params}
	\end{table}
	
\section{Energy-based modeling for the induction motor}
\label{sec:asynchronous}

	We now apply our approach to the Induction Motor (IM). We show that taking the most basic assumptions (sinusoidal and linear motor) we find again the linear model as we did in section \ref{sec:symmetries:linearsinus}.

	\subsection{Deploying the formalism}
		Assuming the squirrel-cage rotor is actually equivalent to three identical wound phases, the generalized coordinates of an IM with three identical stator windings are
		\begin{equation*}
			q = \transpose{(\theta, q_s^a, q_s^b, q_s^c, q_r^a, q_r^b, q_r^c)},
		\end{equation*}
		where $\theta$ is the (electrical) rotor angle and $q_s\abc$ and $q_r\abc$ are the electrical charges in the stator and rotor windings respectively. Their derivatives are
		\begin{equation*}
			\dot{q} = \transpose{(\omega, \is^a, \is^b, \is^c, \ir^a, \ir^b, \ir^c)},
		\end{equation*}
		where $\omega$ is the (electrical) rotor velocity and $\is\abc$ and $\ir\abc$ are the currents in stator and rotor windings respectively. Proceeding as in \ref{sec:hamiltonian:abc}, the generalized momenta are
		\begin{equation*}
			p = \transpose{(\rho, \phis^a, \phis^b, \phis^c, \phir^a, \phir^b, \phir^c)},
		\end{equation*}
		where $\rho$ is the kinetic momentum and $\phi\abc$ and $\phir\abc$ are the flux produced by stator and rotor windings respectively.
		The power exchanges are:
		\begin{itemize}
			\item the electrical power~$\transpose{\us\abc}\is\abc$ provided to the motor by the electrical source, where $\us\abc$ is the vector of voltage drops along the stator winding; this power is associated with the generalized force $\us\abc$
			\item the electrical power~$-\Rs\transpose{\is\abc}\is\abc$ dissipated in the stator resistances $\Rs$; it is associated with the generalized force $-\Rs\is\abc$.
			\item the electrical power~$-\Rr\transpose{\ir\abc}\ir\abc$ dissipated in the rotor resistances $\Rr$; it is associated with the generalized force $-\Rr\ir\abc$.
			\item the mechanical power $-\Tl\frac{\omega}{\np}$ dissipated in the load, where $\Tl$ is the load torque and $\np$ the number of pole pairs; it is associated with the generalized force $-\Tl$.
		\end{itemize}
		
		Using the same method as in \ref{sec:hamiltonian:abc}, we find
		\barCl
			\label{eqn:asynchronous:dynamicHabc}
			\Tderive{\phis\abc} &=& \us\abc - \Rs\is\abc \ase \\
			\Tderive{\phir\abc} &=& -\Rr\ir\abc	\ase \\
			\np\Tderive{\rho} &=& \Te\abc - \Tl, \ase
		\earCl
		where the stator variables are expressed in the stator frame and the rotor variables are expressed in the rotor frame. The current-flux and electro-mechanical relations are also similar,
		\barCl
			\is\abc(\theta, \rho, \phis\abc, \phir\abc) &:=& \Pderive{\HH\abc}{\phis\abc}(\theta, \rho, \phis\abc, \phir\abc) \ans \label{eqn:asynchronous:isabc} \\
			\ir\abc(\theta, \rho, \phis\abc, \phir\abc) &:=& \Pderive{\HH\abc}{\phir\abc}(\theta, \rho, \phis\abc, \phir\abc) \ans\label{eqn:asynchronous:irabc} \\
			\Te\abc(\theta, \rho, \phis\abc, \phir\abc) &:=& -\np\Pderive{\HH\abc}{\theta}(\theta, \rho, \phis\abc, \phir\abc). \ans \label{eqn:asynchronous:Teabc}
		\earCl
		
		Due to the connection scheme of the rotor,
		\begin{equation}
			\ir^a + \ir^b + \ir^c = 0
			\label{eqn:asynchronous:i0}
		\end{equation}
		and the fact that most stators are star-connected (see figure \ref{fig:hamiltonian:star}), it is still interesting to change frame and decouple the $0$\nbdash axis as was done in \ref{sec:hamiltonian:AB}. It is also interesting to express all the variables in the same frame rotating at the synchronous speed $\omega_s$. To do so we define $x_s\dqO := \transpose{\K{\theta_s}}x_s\abc$ and $x_r\dqO := \transpose{\K{\theta_s - \theta}}x_r\abc$ where $\Tderive{\theta_s} := \omega_s$ and
		\begin{equation*}
			\K{\theta} := \sqrt{\frac{2}{3}}\begin{pmatrix}
				\cos{\theta} & \cos{\theta - \frac{2}{3}} & \cos{\theta - \frac{4}{3}} \\
				-\sin{\theta} & -\sin{\theta - \frac{2}{3}} & -\sin{\theta - \frac{4}{3}} \\
				\frac{1}{\sqrt{2}} & \frac{1}{\sqrt{2}} & \frac{1}{\sqrt{2}} \\
			\end{pmatrix}
		\end{equation*}
		Even through the equation will not be preserved, as in \ref{sec:hamiltonian:dq}, we can get similar relations
		\barCl
			\Tderive{\phis\dq} &=& \us\dq - \Rs\is\dq - \J\omega_s\phis\dq \ase \\
			\Tderive{\phir\dq} &=& -\Rr\ir\dq - \J(\omega_s - \omega)\phir\dq \ase \\
			\np\Tderive{\rho} &=& \Te\dq - \Tl \ase
		\earCl
		These are the usual dynamic equations for the IM (see e.g. \cite{Chiasson2005,Krause2002}).
		
		In the $dq$ frame the current-flux and electromechanical relations then read
		\barCl
			\is\dq(\theta, \rho, \phis\dq, \phir\dq) &:=& \Pderive{\HH\dq}{\phis\dq}(\theta, \rho, \phis\dq, \phir\dq) \label{eqn:asynchronous:isdq} \\
			\ir\dq(\theta, \rho, \phis\dq, \phir\dq) &:=& \Pderive{\HH\dq}{\phir\dq}(\theta, \rho, \phis\dq, \phir\dq) \label{eqn:asynchronous:irdq} \\
			\Te\dq(\theta, \rho, \phis\dq, \phir\dq) &:=& -\np\Pderive{\HH\dq}{\theta} + \np\transpose{\ir\dq} \J \phir\dq. \label{eqn:asynchronous:Tedq}
		\earCl
		
	\subsection{Symmetries}
		We now use the motor construction symmetries as in section \ref{sec:symmetries} considering only the case of a sinusoidal induction machine.
		
		So, whatever the angle $\theta$ of the rotor, the energy will be the same, as long as the relative position of the rotor flux space vector with respect to stator flux space vector remains the same. Thus the energy function in the $dq$ frame does not depend on $\theta$.
		
		Rotating the stator and rotor flux space vectors by the same angle $\eta$ preserves the energy, so
		\begin{equation}
			\HH\dq(\rho, \phis\dq, \phir\dq) = \HH\dq(\rho, \R{\eta}\phis\dq, \R{\eta}\phir\dq).
			\label{eqn:asynchronous:invariance}
		\end{equation}
		
		Exchanging two phases on the stator and the rotor and symmetrizing the rotor position also preserves the energy so
		\begin{equation}
			\HH\dq(\rho, \phis\dq, \phir\dq) = \HH\dq(-\rho, \S\phis\dq, \S\phir\dq),
			\label{eqn:asynchronous:symmetry}
		\end{equation}
		with
		\begin{equation*}
			\S := \begin{pmatrix}
				1 & 0 \\
				0 & -1
			\end{pmatrix}.
		\end{equation*}
		
	\subsection{The linear sinusoidal model}
		We consider a second order-polynomial energy function independent on $\theta$ and with magnetic part independent on $\rho$. We keep the simplest expression of the kinetic energy. Such a model is of the form
		\barCl
			\HH_l\dq &:=& \frac{1}{2\Jl\np^2}\rho^2 + a + b\phis\dq + c\phir\dq \ane \\
			&+& \transpose{\phis\dq}D\phis\dq + \transpose{\phis\dq}E\phir\dq + \transpose{\phir\dq}F\phir\dq,
		\earCl
		where $a \in \mathbb{R}$, $(b, c) \in (\mathbb{R}^2)^2$ and $(D, E, F) \in (\mathcal{M}_2(\mathbb{R}))^3$.
		
		The equation \eqref{eqn:asynchronous:invariance} implies that $b = c = (0,0)$ and $D$, $E$ and $F$ commute with the rotations. So
		$(D, E, F) \in \left\{\alpha \I + \beta \J, (\alpha, \beta)\in \mathbb{R}^2\right\}$ where $\I \in \mathcal{M}_2(\mathbb{R})$ is the identity matrix and $\J$ was defined in \ref{sec:hamiltonian:dq}. Due to \eqref{eqn:asynchronous:symmetry} $D$, $E$ and $F$ are colinear with $\I$ because $\J$ does not commute with $\S$, hence the energy function is of the form
		\begin{equation}
			\HH_l\dq := \frac{1}{2\Jl\np^2}\rho^2 + a + d\transpose{\phis\dq}\phis\dq + e\transpose{\phis\dq}\phir\dq + f\transpose{\phir\dq}\phir\dq.
		\end{equation}
		We can choose freely $a = 0$ as the energy function is defined up to a constant. We define $\sigma$, $\Lm$, $\Ls$ and $\Lr$ by the implicit relations (it can be checked that it is invertible when it is defined)
		\begin{equation*}
			\Lr\Ls\sigma = \Ls\Lr - \Lm^2
		\end{equation*}
		\begin{equation*}
			d = \frac{1}{2\Ls\sigma} \qquad
			e = -\frac{2\Lm}{2\Lr\Ls\sigma} \qquad
			f = \frac{1}{2\Lr\sigma}
		\end{equation*}
		Thus, the energy function reads
		\barCl
			\HH\dq &:=& \frac{1}{2\Jl\np^2}\rho^2 + \frac{\Lm}{2\Ls\Lr\sigma}\transpose{(\phis\dq - \phir\dq)}(\phis\dq - \phir\dq) \ane \\
			&+& \frac{\Lr-\Lm}{2\Ls\Lr\sigma}\transpose{\phis\dq}\phis\dq
			+ \frac{\Ls-\Lm}{2\Ls\Lr\sigma}\transpose{\phir\dq}\phir\dq.
		\earCl
		Applying \eqref{eqn:asynchronous:isdq} and \eqref{eqn:asynchronous:irdq} one gets the current-flux relations
		\barCl
			\Ls\Lr\sigma\is\dq = \Lm(\phis\dq - \phir\dq) + (\Lr - \Lm)\phis\dq \ane \\
			\Ls\Lr\sigma\ir\dq = \Lm(\phir\dq - \phis\dq) + (\Ls - \Lm)\phir\dq. \ane
		\earCl
		Inverting these equations and taking into account the electro-mechanical torque is $\Te = n\transpose{\ir\dq}\J\phis\dq$, the usual relations (see e.g. \cite{Chiasson2005,Krause2002}) are easily identified. Therefore we recovered the linear sinusoidal model for the IM without the tedious microscopic approach.


\bibliographystyle{phmIEEEtran}
\bibliography{biblio}

%
%
%
%

\end{document}